\documentclass[10pt,a4paper]{article}
  \usepackage{amssymb,amsmath,exscale}
  \usepackage{bm}
  \newcommand{\bd}{ \pmb d }
 
    \usepackage{bbm}
 
\newcommand{\intzt}{\int_0^\zt}
\newcommand{\ZA}{\mathcal{A}_\zt}
\newcommand{\ZB}{\mathcal{B}_\zt}
\newcommand{\ZC}{\mathcal{C}_\zt}
 \newcommand{\PP}{\mathcal{P}}
 
 \newcommand{\ZNN}{\mathcal{N}}
\newcommand{\dd}{d_1}
\newcommand{\D}{d_2}
\newcommand{\upzt}{{u^+_\zt}}
\newcommand{\wpzt}{{w^+_\zt}}
\newcommand{\intozt}{\int_0^{\zt}}
\newcommand{\intztt}{\int_\zt^t}
\newcommand{\uoptzt}{u^+_\zt}
\newcommand{\woptzt}{w^+_\zt}
\newcommand{\heaviside}{\mathbbm{1}}
\newcommand{\ztu}{{t}}
 
\usepackage[utf8]{inputenc} 
\usepackage{makeidx}
  \usepackage{bbm}
  \usepackage{bm}
\usepackage{amsmath} 
\usepackage{amsfonts}
\usepackage{amssymb}
   \usepackage{amssymb,amsmath,exscale}
 
\usepackage{epsfig}
\usepackage{graphics}

  \input definit.tix
 
\author{
L. Pandolfi\thanks{Dipartimento di Scienze Matematiche ``Giuseppe Luigi Lagrange'', Politecnico di Torino, Corso Duca degli Abruzzi 24, 10129 Torino, Italy (retired) (luciano.pandolfi@formerfaculty.polito.it)}
}
\title{The quadratic tracking problem for systems with persistent memory in $\zzr^d$}

\begin{document}

\maketitle

{\bf Abstract:}   The classical   quadratic regulator problem has rarely been studied for systems with persistent memory until recent times. In this paper we study the quadratic tracking problem on a \emph{ finite time horizon} for a system described by a controlled linear Volterra integrodifferential equation in $\zzr^d$. We use the Fredholm equation approach and we derive the synthesis of the optimal control  in terms of a Riccati differential equation, independent of the reference signal, and of two equations which instead depend on the reference signal.  In the final section we introduce a representation of the system in a state space   \emph{of finite memory  } (as the tracking problem under study) and  we show that the equations used in the synthesis  can be formulated as a differential system in this space. This fact has to  be contrasted with the    semigroup approach  which requires that the system be recasted in a space of infinite memory.

\medskip

{\bf key words:} Systems with persistent memory, Volterra integrodifferential equations, tracking problem

\section{Introduction} 	 
 
We study the tracking problem for a system with persistent memory, described as follows. The system is
\begin{subequations}
\begin{multline}
\ZLA{eq:SisteVolte0}
w'(t)=Aw(t)+\intt N(t-s) w(s)\ZD s+ Bu(t)\,,\qquad t\geq 0\,,\\
 w(0)=\hat\xi\in \zzr^d  
\end{multline}
where $A $ is a constant $d\times d $ matrix and 
$N(t)$ is a continuous $d\times d$ matrix.
 
 The $d\times m$  \emph{input matrix}     $B$ is constant     and the \emph{control} $u$ is a     square integrable $\zzr^m$-valued function.   
 
We assign a  \emph{continuous reference signal} $y(t)\in \zzr^p$. The \emph{tracking problem} is the study of the optimal control $u^+$ defined by
\begin{align}
\nonumber
u^+&={\rm arg \min}\, J_0 (\hat\xi;u)\,,\\
\ZLA{eq:SisteCOSTO0}
&J_0 (\hat\xi;u)=\int_0^T \left [\, \|Cw(t)-y(t)\|^2 +\|u(t)\|^2\,\right ]\ZD t \,.
\end{align}
\end{subequations}
The final time $T$ is fixed, the $p\times d$ \emph{output matrix}  $C$ is constant and $w$ solves~(\ref{eq:SisteVolte0}). The norms   in the cost $J_0$ are those, respectively, of $\zzr^p$ and of $\zzr^m$. This is not explicitly indicated without any risk of confusion.

References to this classical problem  are discussed in Sect.~\ref{sect:REFERENCES}. Here we recall a fact that is well understood: it is convenient to introduce a family of optimal control problems  parametrized by $\zt$, any initial time in $[0,T)$. So, we study the family of the optimization problems
 \begin{subequations}
\begin{align}
\nonumber
w'(t)&=Aw(t)+\int_\zt ^t N(t-s) w(s)\ZD s+ Bu(t)\\
\ZLA{eq:SisteVoltezt}
&+\int_0^\zt N(t-s)\tilde \xi(s)\ZD s\,,\qquad w(\zt)=\hat\xi \in \zzr^d\,, 
 \\
\nonumber
u&={\rm arg \min}\, J_\zt ((\hat \xi;\tilde \xi\ZCD),u)\,,\\ 
\ZLA{eq:SisteCOSTOzt}
&J_\zt ((\hat \xi;\tilde \xi\ZCD),u)=\int_\zt^T \left [\, \|Cw(t)-y(t)\|^2 +\|u(t)\|^2\,\right ]\ZD t\,.
\end{align}
\end{subequations}
The initial condition at the time $\zt\in [0,T)$ of the equation~(\ref{eq:SisteVoltezt}) is $\hat \xi$ when $\zt=0$ while it is the pair $(\hat \xi,\tilde \xi\ZCD)$ when $\zt\in (0,T)$. The function $\tilde \xi$ needs not be the restriction of the solution of~(\ref{eq:SisteVolte0}) to $[0,\zt]$. It is an \emph{arbitrary} square integrable function. It is convenient to introduce the notation $\Xi_\zt$:
\[\left\{\begin{array}{ll}
\mbox{if $\zt=0$ then}&
\Xi_0=\hat \xi\in \zzr^d \\
\mbox{if $\zt\in (0,T)$ then}&\Xi_\zt=(\hat \xi,\tilde \xi\ZCD)\in M^2_\zt= \zzr^d\times L^2(0,\zt;\zzr^d)\,.
\end{array}\right.
\]

The existence of a unique optimal control is clear 
 and the standard variational method provides a representation for the optimal control (see Sect~\ref{SecVARIATIONAL}). Our goal is to derive a system of equations (including a suitable version of the Riccati equation) from which the optimal control can be directly computed.
 
\subsection{\ZLA{sect:REFERENCES}Organization of the paper and comments on the literature}
 
Systems with persistent memory, described as in~(\ref{eq:SisteVolte0}) and~(\ref{eq:SisteVoltezt}) are often encountered in applications even when   $w$ belongs to an infinite dimensional Hilbert spaces (for example, as models of   viscoelastic bodies or diffusion processes, see~\cite{PandolfiLIBRO21,PruessLIBRO1993}) but in this paper we confine ourselves to study the   case~$w\in\zzr^d$.

In most of the applications the matrices $A$ and $N(t)$ commute, and even more, quite often   $N(t)=n(t)A$ where $n(t)$ is a real valued function, but we don't use this assumption in this paper: \emph{the matices need not commute.}

The paper is organized as follows. In Sect.~\ref{SecVARIATIONAL} we derive a convenient representation of the solutions of Eq.~(\ref{eq:SisteVoltezt}). Then we use the standard variational approach to derive a formula for the optimal control of problem~(\ref{eq:SisteVoltezt})-(\ref{eq:SisteCOSTOzt}).  The focus in this section is on the study of  the regularity of the optimal control. In 
section~\ref{sectTheDI} we show a suitable version of the dissipation inequality (DI) (both in integral and differential form) which holds for the problem~(\ref{eq:SisteVoltezt})-(\ref{eq:SisteCOSTOzt}) and the corresponding equality  satisfied by the optimal control.
In Sect.~\ref{sec:FinalRepre}
we couple with the regularity properties derived in Sect.~\ref{SecVARIATIONAL} and we derive a systems of equations which allows complete identification of the optimal control in terms of a feedback control and a term which solely depend on the reference signal $y$. 

In the final section~\ref{OpeFORMequa} we show a state space representation 
of Eq.~(\ref{eq:SisteVolte0}) (and so also of~(\ref{eq:SisteVoltezt})) in a memory space with \emph{finite memory.}  State representations of Eq.~(\ref{eq:SisteVolte0}) have been proposed by using 
semigroup theory (see~\cite{Amendola2012,ENGELnagelLIBOsemigruppi,PandolfiLIBRO21}). The semigroup approach require that the memory of the system be enlarged to the entire past, up to $t=-\ZIN$. So, the semigroup  approach is not natural when studying optimization problems on a finite time horizon. For this reason in section~\ref{OpeFORMequa}  we introduce a 
  state space representation of the system with memory in a space with finite memory and we recast in this space the synthesis of the optimal control derived in Sect.~\ref{sec:FinalRepre}.  

As documented in~\cite{PandolfiLIBRO21},  Systems~(\ref{eq:SisteVolte0}) and~(\ref{eq:SisteVoltezt})  have been mostly studied from the point of view of stability and controllability. Preliminary results on the   quadratic regulator problem are   in~\cite{pritchardLUO96}. After this paper, the study of the quadratic regulator was not pursued for a long time and only recently  people in probability and game theory found interest in this problem (see for example~\cite{ShuoLin23,WANG18ESAIm} 
and references therein). In the special case $A=0$ and $y=0$, the paper~\cite{PandolfiVOLTERRAieee18} derived a characterization of the optimal control based on a suitable version of   the Riccati equation The results in~\cite{PandolfiVOLTERRAieee18} have been extended to a class of distributed systems  (with $y=0$ and scalar memory kernel) in~\cite{acquistaBUCCI24}. 

A classical presentation of the tracking problem (for systems without memory) is~\cite{AthansFALB1966}. A selfcontained extension to a class of (memoryless) distributed systems is in~\cite{PandolfiPRIOLAtracking}.

In the quadratic tracking problem, the reference signal $y$ is known in advance (as we assume in the present paper). We mention that in certain applications instead an exterior signal  has to be recursively identified and at the same time it is required either that it is tracked or, when it is a disturbance, it is required that its influence on the system is reduced as much as possible (this is the problem of disturbance rejection).   These goals are achieved either by using Kalman   and Luemberger filters (see~\cite{DAVISfoundatCONTROL02,SontagLIBRO98}) or via a deterministic  approach to the identification of an exterior  signal    proposed by Krasovski (the book~\cite{MaksimovLIBRO02} extends Krasovski method also to a class of (memoryless) distributed systems). See~\cite{FagnaniPandJIIPP03,MFAGnanIPandoLfiIPvoltyeDecon03} for extensions to systems with memory and~\cite{FagnaniMaksPandolfi,PandolfiINTjCONTreject} for an application to  disturbance rejection (for finite dimensional, memoryless systems).

\section{\ZLA{SecVARIATIONAL}Representation of the solutions and the regularity of the optimal control}
We denote $Z(t)$ the unique $d\times d $ matrix which solves
\begin{equation}\ZLA{Eq:dellaFUNZzMaiusc}
Z'(t)=A^* Z(t)+\intt N^*(t-s)Z(s)\ZD s\,,\qquad Z(0)=I\ \mbox{(the identity matrix)}\,.
\end{equation}
It is easy to compute that for  $r\in (0,t)$  the solution $w$ of~(\ref{eq:SisteVoltezt}) satisfies
\begin{multline*}
\dfrac{\ZD}{\ZD r} Z^*(t-r)w(r)=Z^*(t-r) \int_\zt^r N(r-s)w(s)\ZD s\\
-\left [  \int_0^{t-r} Z^*(s)N(t-r-s)\ZD s   \right ]w(r)+Z^*(t-r)\left [Bu(r)+\int_0^\zt N(r-s)\tilde \xi(s)\ZD s\right ]
\end{multline*}
so that
\begin{multline}
\ZLA{eq:FormVARconst}
w(t)=w(t;\zt,\Xi_\zt,u)=\left [Z^*(t-\zt)\hat\xi+\intztt Z^*(t-r)\intozt N(r-s)\tilde\xi(s)\ZD s\,\ZD r\right ]
\\
+\intztt Z^*(t-r)Bu(r)\ZD r\qquad t\in[\zt,T]\,.
\end{multline}
This is the variation of constants  formula for the equation~(\ref{eq:SisteVoltezt}).  

We extend the definition of $w(t)=w(t;\Xi_\zt,u)$ to $[0,\zt]$ by putting
\begin{equation}
\ZLA{eq:DefiESTEdiW}
w(t;\zt, \Xi_\zt,u)=\left\{\begin{array}
{l}
\mbox{the function~(\ref{eq:FormVARconst})   if $t\in [\zt,T]$ } \\
\mbox{$\tilde \xi$ if $t\in (0,\zt)$}\,.
\end{array}\right.
\end{equation}

Note that $w(\cdot;\zt, \Xi_\zt,u)\in C([\zt,T],\zzr^d)$  and  $w(\zt;\zt,\Xi_\zt,u)=\hat\xi$.

We need the following obvious consequence of~(\ref{eq:FormVARconst}):
\begin{Lemma}\ZLA{eqDERIVrispeZR}
Let $\xi\in C([0,T];\zzr^d)$ and let $\tilde\xi =\xi_{|_{(0,\zt)}}$, $\hat\xi=\xi(\zt)$.  If $u$ is continuous then $(t,\zt)\mapsto w(t;\zt,\Xi_\zt,u)$ is of class $C^1$ on the angle $t\geq \zt\geq 0$ (the derivatives admits finite limits for $(t,\zt)\to (s,s)$  from inside the angle). 
\end{Lemma}

The representation~(\ref{eq:FormVARconst}) shows that, for every $T>\zt$, the transformation $(\Xi_\zt,u)\mapsto w$ is linear and continuous from $M^2_\zt\times L^2(\zt,T;\zzr^d)$ to 
$C([\zt,T],\zzr^d)$ and so the minimization problem~(\ref{eq:SisteVoltezt})-(\ref{eq:SisteCOSTOzt}) admits a unique optimal control.
 
The optimal control is denoted $\uoptzt $ and the corresponding solution of~(\ref{eq:SisteVoltezt}) is denoted $\woptzt$. These function do depend on the initial datum $\Xi_\zt$ and when needed for clarity they are denoted 
$\uoptzt (\cdot;\Xi_\zt)$ and $\woptzt (\cdot;\Xi_\zt)$. We use the standard variational method, i.e. we insert~(\ref{eq:FormVARconst}) in the cost and we put equal zero the Gateaux derivative (respect to $u$). We find that the optimal control is characterized by
\begin{equation}
\ZLA{eq:expreOptCONTviaCoistato}
\uoptzt (t;\Xi_\zt)= -B^*p_\zt(t)\,,\quad p_\zt(t)=  \int_t^T Z(s-t)C^*\left [C \woptzt  (s;\Xi_\zt)-y(s)\right ]\ZD s \,.
\end{equation}
We recall $0\leq \zt\leq t\leq T$ and we note that $p_\zt$ does depend also on the initial condition $\Xi_\zt$.

It is seen that $p=p_\zt$ satisfies
\begin{equation}\ZLA{eq:delCOSTATO}
p'(t)=-A^*p(t)-\int_t^T N^*(s-t)p(s)\ZD s-C^*\left [C \woptzt  (t;\Xi_\zt)-y(t)\right ] 
\end{equation}
with the final condition $p(T)=0$.
  As we are not going to use the differential equation~(\ref{eq:delCOSTATO}), we leave its derivation   to the interested reader.

 \begin{Remark}\ZLA{remaSect1EFFETTOpassato}{\rm
Let Eq.~(\ref{eq:SisteVoltezt}) contains an additive term $F(t)$. Then the additive term
\[
\int_0^t Z^*(t-s)F(s)\ZD s
\]
appears in~(\ref{eq:FormVARconst}). The additive term $F(t)$ has the sole effect os changing $y(t)$ in~(\ref{eq:expreOptCONTviaCoistato}) to 
\[
y_1(t)=C\intt Z^*(t-s)F(s)\ZD s-y(t)\,.
\]
So, the contribution if $F$ is absorbed by the reference signal and the presence of $F$ does not lead to a more general problem. 

The interest of this obvious observation stems from the fact when $F(t)=\int_{-\ZIN}^0 N(t-s)\xi(s)\ZD s$ the additive term $F$ represents the memory of the configuration of the system prior to the time  $0$.\zdia
 }\end{Remark}
 
 The first consequence of the representation~(\ref{eq:expreOptCONTviaCoistato}) of the optimal control is as follows:
 \begin{Lemma}\ZLA{Lemma:ContiDIFFEoptimalCONTROL}
 For every fixed $\zt\in [0,T)$ the function $t\mapsto \uoptzt (t;\Xi_\zt)$ is continuously differentiable (in particular, it is continuous) on the interval $[\zt,T]$ (with finite directional derivative at the end points). 

 \end{Lemma}

In order to simplify the notation, we rename $ Y_{0,\zt}$ the bracket in~(\ref{eq:FormVARconst}):
\[
Y_{0,\zt}(t;\zt)=Z^*(t-\zt)\hat\xi+\intztt Z^*(t-r)\intozt N(r-s)\tilde\xi(s)\ZD s\,\ZD r
\]
($Y_{0,\zt}$ depends   on $\Xi_\zt$).

We replace $u$ with $-B^*p_\zt$ in~(\ref{eq:FormVARconst}) so that the solution      is $\wpzt$. Then we replace~(\ref{eq:FormVARconst})   in~(\ref{eq:expreOptCONTviaCoistato}). We get
 
\begin{subequations}
\begin{eqnarray}\ZLA{FredhEQUAdip}\nonumber
p_\zt (t)&=& \int_t^T Z(s-t)C^*\left [CY_{0,\zt}(\zt)-y(s)\right ]\ZD s
\\
\ZLA{FredhEQUAdipA}&\quad&- \int_t^T Z(s-t)C^*\left [C \int_\zt^sZ^*(s-r)B\left (B^*p_\zt(r)\right )\ZD r\right  ]\ZD s\\
\nonumber
&=& \int_t^T Z(s-t)C^*\left [CY_{0,\zt}(s;\zt)-y(s)\right ]\ZD s\\
\nonumber
&\quad & -\int_\zt^t\left [\int_t^T Z(s-t)C^*CZ^*(s-r)B\ZD s\right ] \left (B^* p_\zt(r)\right )\ZD r\\
\ZLA{FredhEQUAdip}
&\quad & -\int_t^T\left [\int_r^T Z(s-t)C^*CZ^*(s-r)B\ZD s\right ] \left (B^* p_\zt(r)\right )\ZD r
\,.
\end{eqnarray}
 
\end{subequations}
 
We introduce the function $Y_\zt(t)$ (which depends also on   $\Xi_\zt$):
\begin{equation}\ZLA{eq:LaFORMperY}
Y_\zt(t)= \int_t^T Z(s-t)C^*\left [CY_{0,\zt}(s;\zt)-y(s)\right ]\ZD s\,.
\end{equation}
 Eq.~(\ref{FredhEQUAdip}) is a Fredholm integral equation of the second kind for $p_\zt$:

\begin{equation}
\ZLA{FormuladiFREDhPERp}
p_\zt(t)=Y_\zt(t)-\int_\zt^T K(t,r)\left (B^*p_\zt(r)\right )\ZD r
\end{equation}
where
 
\begin{align*}
K(t,r)&=\tilde K(t,r)B= \left \{\heaviside(t-r)\left [\int_t^T Z(s-t)C^*CZ^*(s-r) \ZD s\right ]\right .\\
&\left.+\heaviside(r-t)\left [\int_r^T Z(s-t)C^*CZ^*(s-r) \ZD s\right ]\right \}B\,,\quad (t,r)\in[\zt,T]\times [\zt,T]\,.
 \end{align*}
For most of clarity we specify that $\tilde K(t,r)$ (the brace) is a $d\times d$ matrix while $K(t,r) \in\zzr^{m\times d}$.
Furthermore we note 
\begin{equation}
 \ZLA{EqCONDIfinK}
 K(T,r)=0\  \forall  r\,,\qquad K(t,T)=0\ \forall t\,.
 \end{equation} 
 
Thanks to the fact that $Z(t)\in C^1([0,+\ZIN))$ (with finite directional derivative at $0$), we see that the kernel $K(t,r)$ is continuous on $[\zt,T]\times[\zt,T]$, continuously differentiable if $t\neq r$ and with bounded derivatives both in $t>r>0$ and in $0<r<t$. And so $K(t,r)$ is even Lipschitz 
continuous.
 
 We have:
\begin{Lemma}\ZLA{LemmaESIuniSOLUfredP}
Let $\zt\in[0,T]$ be fixed.
The Fredholm integral equation~(\ref{FormuladiFREDhPERp}) admits a unique solution for every $Y_\zt\in L^2(0,\zt)$ (note that the function $Y_\zt$ in~(\ref{eq:LaFORMperY}) is even continuous).
\end{Lemma}
\zProof First we prove that the solution is unique. Linearity of the problem shows that it is sufficient to prove that when $Y_\zt=0$ the equation admits the unique solution $p=0$.

If $Y_\zt=0$ and $B^*p_\zt=0$ then equality~(\ref{FormuladiFREDhPERp}) shows that $p_\zt=0$. So, it is sufficient to prove that the following equation admits the sole solution $v=0$:
\[
v(t)+\left (L v\right )(t)=0\,,\qquad \left (Lv\right )(t)= \int_t^T B^* Z(s-t)C^*\left [  \int_\zt^sCZ^*(s-r)Bv(r)\ZD r\right  ]\ZD s\,.
\]
This follows since   the operator $L$ is selfadjoint positive in $L^2(\zt,T;\zzr^m)$.

For this same reason, the operator $I+L$ from $L ^2(\zt,T;\zzr^m)$ to itself is surjective, so that the solution $B^*p$ (and so also $p$) exists for every square integrable function $Y$.\zdia

The theory of the Fredholm integral equation asserts that the unique solution has the following representation:
\begin{equation}
\ZLA{eq:FredhRISOLTA}
p_\zt(t)=Y_\zt (t)-\int_\zt^T R(t,r;\zt) Y_\zt(r)\ZD r\,.
\end{equation}
The \emph{resolvent kernel} $R(t,s;\zt)\in\zzr^{d\times d}$ is the unique solution of
 
\begin{subequations}
\begin{align}
\ZLA{eq:delRIsolventeA}
R(t,r;\zt)&=-\int_\zt^T K(t,\nu)B^*R(\nu,r;\zt)\ZD\nu+K(t,r)B^*\\
\nonumber
&=
- \int_t^T Z(s-t)C^*\left [C \int_\zt^sZ^*(s-\nu )B\left (B^*R(\nu,r;\zt) \right )\ZD \nu\right  ]\ZD s\\
 \ZLA{eq:delRIsolventeB} 
&\hskip  5cm +K(t,r)B^*\,,\qquad t\in[\zt,T]
\end{align}
\end{subequations}
(see~(\ref{FredhEQUAdipA}) and~(\ref{FredhEQUAdip}) and recall $K(t,r)B^*=\tilde K(t,r)BB^*\in\zzr^{d\times d}$).
In this equation, both $t$ and $r$ belong to $ [\zt,T]$ and
\begin{equation}
 \ZLA{EqCONDIfinR}
 R(T,r;\zt)=0 \,,\qquad R(t,T;\zt)=0\,,\qquad R(t,r;T)=0
 \end{equation}
 (the first two equality follows from~(\ref{EqCONDIfinK}) while the last equality since $\zt\leq t\,,\ r\leq T$).

\subsection{The regularity of the resolvent kernel $R(t,r;\zt)$ and its consequences}

 A known fact, from the Theory of the Fredholm integral equation, is that $(t,r)\mapsto R(t,r;\zt)$ is continuous on $[\zt,T]\times [\zt,T]$, as $K(t,r)$. We need a more precise statement.
 
 \emph{In the following statements, when stating that a function is ``of class $C^1$ on the closure of a region'' we intend that the derivatives have continuous extensions to the boundary.}
In the following proof we repeatedly use the following obvious observation:
 
\begin{Lemma}\ZLA{lemmaCONTIinteg}
Let $\Phi(x,s,y)$ be a bounded function of   three real variables and let $x$ and  $s$ belong to the same interval $[a,b]$  while $y\in [c,d]$ (possibly a different interval). We assume  
1) the function $(x,s,y)\mapsto \Phi(x,s,y)$ is continuous on the set $\{(x,s,y)\,\ x\neq s\}$.; 2) the directional limits of   $x\to \Phi (x,s,y)$ exist and are finite   when $x$ tends to   to $s$ (possibly different directional limits).

The function 
\[
\Psi(x,y)=\int _a^b \Phi(x,s,y)\ZD s \int _a^x \Phi(x,s,y)\ZD s+\int_x^b  \Phi(x,s,y)\ZD s
\]
is continuous.
 
\end{Lemma}

\begin{Theorem}\ZLA{TheoremPROPdiR}The following properties hold:
\begin{enumerate}
\item\ZLA{I1TheoremPROPdiR} There exists a number $M_T$, which does not depend on $\zt\in [0,T]$, such that
\[
\|R(t,r;\zt)\|<M_T\qquad\zt\in[0,T]\,,\ (t,r)\in[\zt,T]\times[\zt,T] 
\]
(the norm is that of  the $d\times d$ matrices).
\item\ZLA{I2TheoremPROPdiR}
The function $(t,r)\mapsto R(t,r;\zt)$ has the following properties:
\begin{enumerate}
\item\ZLA{I21TheoremPROPdiR} it is 
Lipschitz continuous on $[\zt,T]\times [\zt,T]$  uniformly in $\zt$, i.e. there exists $M=M_T$ such that
\begin{equation}
\ZLA{eq:LipContRunizt}
|R(t,r;\zt)-R(t_1,r_1;\zt)|<M_T\left \{\, |t-t_1|+|r-r_1|\,\right \}  \,.
\end{equation}
\item\ZLA{I22TheoremPROPdiR} 
it is of class $C^1$ on the set $\{(t,r)\in[\zt,T]\times [\zt,T]\,, \ t\neq r\}$ and the derivatives admits limits for $( t,r)\to (s,s)$ either from above or from below the bisectrix.
\end{enumerate}
\item\ZLA{I3TheoremPROPdiR} Let $t_0$ and $r_0$ be such that $0<t_0\leq T$, $0<r_0\leq T$ and let $\zt_M=\min\{t_0,\,r_0\}$. The function
\[
  R(t_0,r_0;\cdot) \in C^1([0,\zt_M])\,. 
\]
 
\end{enumerate}
\end{Theorem}
\zProof \emph{We prove the statement~\ref{I1TheoremPROPdiR}. }  We fix $\zt\in[0,T]$. 
We introduce the operator $\Lambda_\zt $:
\[
(\Lambda_\zt v)(t)=C\int_\zt^tZ^*(t-\nu) Bv(\nu)\ZD \nu\,,\quad \Lambda_T\in\mathcal{L}(L^2(\zt,T;\zzr^m),L^2(\zt,T;\zzr^p))\,.
\]
Let  $r\in[\zt,T]$. The equality~(\ref{eq:delRIsolventeB})
 gives
 \[
B^*R(\cdot,r;\zt)+\left (\Lambda_\zt^*\Lambda_\zt\right )\left ( B^*R(\cdot,r;\zt)\right )(\cdot)=B^*K(\cdot,r)B^* =B^*\tilde K(\cdot,r)BB^*\,.
 \]
The function $K(t,r)$ does not depend on $\zt$ and it is bounded. It follows  the existence of $M_0$  such that
 \[
\| B^*R(\cdot,r;\zt)\|_{L^2(\zt,T;\zzr^d)}<M_0 \quad \forall \zt\in [0,T]\,,\ r\in[\zt,T]
\,.
 \]
So,  the equality~(\ref{eq:delRIsolventeA}) gives
\[
\|R(t,r;\zt)\|\leq \sup_{(t,r)\in[\zt,T]^2}\|K(t,r)B^*\|+M_0\sup_{t\in[\zt,T]} \left [\int_\zt^T\|K(t,\nu)\|^2\ZD\nu\right ]^{1/2}=M\,.
\]
The number $M$ does not depend either on $\zt\in[0,T]$ or on $(t,r)\in[\zt,T]\times[\zt,T]$ as wanted. 

\emph{We prove the statement~\ref{I2TheoremPROPdiR}.} We fix $t_0$, $r_0$ and $\zt_0$ and we consider the incremental quotient   respect to the variable $t$ at $(t_0,r_0)$:
\begin{multline*}
\dfrac{\Delta(t_0,r_0; h)}{h}=\dfrac{R(t_0+h,r_0;\zt_0)-R(t_0,r_0,\zt_0 )}{h} \\
=-\int_{\zt_0}^T \dfrac{K(t_0+h,\nu )-K(t_0,\nu )}{h}B^*R(\nu,r_0;\zt_0)\ZD\nu\\
+
\dfrac{K(t_0+h,r_0 )-K(t_0,r_0 )}{h}B^*\,.
\end{multline*}

Lipschitz continuity of $K$ and the boundedness property in the statement~\ref{I1TheoremPROPdiR} shows the existence of $M(T)$ such that
\begin{equation}\ZLA{eq:LipContRuniztT0}
\left |\, R(t_0+h,r_0;\zt_0)-R(t_0,r_0;\zt_0)\, \right |<M(T) h 
\end{equation}
(this is inequality~(\ref{eq:LipContRunizt}) in the case $r=r_1$).
Then we use
\[
\lim _{h\to 0^+} \dfrac{K(t_0+h,\nu)-K(t_0,\nu)}{h} =K_t(t_0,\nu) \quad\mbox{$ {\rm a.e}\  t\in [\zt_0,T]$ (i.e. if $\nu\neq t_0$)}
\]
and boundedness of the incremental quotient:   we can exchange limit and integral and we get
\[
R_t(t_0,r_0;\zt_0)=K_t(t_0,r_0 )-\int _{\zt_0}^T K_t(t_0,\nu) B^* R(\nu,r_0;\zt_0)\ZD\nu\,.
\]
We apply Lemma~\ref{lemmaCONTIinteg} to the function $\Phi(t_0,r_0,\nu)=K_t(t_0,\nu)B^*R(\nu,r_0;\zt_0)$ (i.e. with $x=t_0$, $s=\nu$, $y=r_0$).
We see  that $(t_0,r_0)\mapsto R_t(t_0,r_0;\zt_0)$ is a continuous function  when $t_0\neq r_0$. 

Now we consider the incremental quotient respect to $r$. We consider it as a function of $t$:

\begin{align*}
\dfrac{\Delta(t, r_0;h)}{h}=\dfrac{R(t ,r_0+h;\zt_0)-R(t ,r_0;\zt_0)}{h}\\
=-\int _{\zt_0}^T K(t_0,\nu)B^*\dfrac{\Delta(\nu,r_0;h)}{h}\ZD\nu+\dfrac{K(t ,r_0+h)-K(t_0,r_0)}{h}B^*\,.
\end{align*}
We solve this Fredholm integral equation on $[\zt_0,T]$ respect to the unknown $(1/h)\Delta(\cdot,r_0;h)$:
\[
\dfrac{\Delta(t ,r_0;h)}{h}=\dfrac{K(t ,r_0+h)-K(t ,r_0)}{h}B^*-
\int _{\zt_0}^T R(t ,\nu;\zt_0)
\dfrac{K(\nu,r_0+h)-K(\nu,r_0)}{h}B^*\ZD\nu\,.
\]
This equality holds for every $(t ,r_0)\in [\zt_0,T]\times [\zt_0,T]$ and, as above, it implies the existence of a constant $M(T)$ such that
\begin{equation}\ZLA{eq:LipContRuniztR0}
|R(t ,r_0+h;\zt_0)-R(t ,r_0;\zt_0)|<M(T)h\,.
\end{equation}
It implies also the equality
\[
R_r(t,r_0;\zt_0)=K_r(t,r_0)B^*-\int _{\zt_0}^T R(t,\nu;\zt_0)K_r(\nu,r_0)\ZD \nu\,.
\]
This derivative too is a continuous function of $(t_0,r_0)$ for $t_0\neq r_0$.

The uniform Lipschitz continuity~(\ref{eq:LipContRunizt}) follows by combining~(\ref{eq:LipContRuniztT0}) and~(\ref{eq:LipContRuniztR0}).

\emph{We prove the statement~\ref{I3TheoremPROPdiR}.} In this part of the proof we use the boundedness  property proved in the statement~\ref{I1TheoremPROPdiR} and
  the Lipschitz continuity property~(\ref{eq:LipContRunizt}).
 
 We fix $\zt_0 \in [0,T)$ and we prove continuity at $\zt_0$ (the right continuity if $\zt_0=0$).  We consider the difference
\begin{multline}\ZLA{eqA:I3TheoremPROPdiR}
\Delta(t ,r ;h)=R(t,r;\zt_0+h)-R(t,r;\zt_0)\\
=-\int _{\zt_0+h}^T K(t,\nu)B^*R(\nu,r;\zt_0+h)\ZD\nu+\int _{\zt_0}^T K(t,\nu)B^* R(\nu,r;\zt_0)\ZD\nu\,.
\end{multline}
In order that this equality make  sense it must be $\max\{\zt_0,\, \zt_0+h\}\leq \zt_M=\min\{t,\,r\}$.

We use~(\ref{eqA:I3TheoremPROPdiR}) and   we prove continuity first. Then we prove differentiability. In the proof of differentiability   the Lipschitz continuity property~(\ref{eq:LipContRunizt}) is used.

\begin{description}
\item[\bf Continuity:] First we prove left continuity (when $\zt_0>0$). Let $h=-\ZSI<0$. We rewrite~(\ref{eqA:I3TheoremPROPdiR})  as the following Fredholm integral equation in the variable $t\in [\zt_0,T]$:
\begin{multline*}
\Delta(t,r;-\ZSI)=
R(t,r;\zt_0-\ZSI)-R(t,r;\zt_0)\\
=
-\int _{\zt-\ZSI}^{\zt_0} K(t,\nu)B^* R(\nu,r;\zt_0-\ZSI)\ZD\nu
-\int_{\zt_0}^T K(t,\nu)B^*\Delta(\nu,r;-\ZSI)\ZD\nu
\end{multline*}
so that

\begin{multline}
\ZLA{eqB:I3TheoremPROPdiR}
\Delta(t,r;-\ZSI)=-\int _{\zt_0-\ZSI}^{\zt_0} K(t,\nu)B^* R(\nu,r;\zt_0-\ZSI)\ZD\nu\\
+\int _{\zt_0} ^T R(t,s;\zt_0)\left [\int _{\zt_0-\ZSI} ^{\zt_0} K(s,\nu)B^* R(\nu,r;\zt_0-\ZSI)\ZD\nu\right ]\ZD s\,.
\end{multline}
The boundedness property in the Statement~\ref{I1TheoremPROPdiR} shows that
\[
\lim _{\ZSI\to 0^+}\Delta(t,r;-\ZSI)=0
\]
and the limit is uniform on the set $\{(t,r)\,:\ t\geq \zt_0, \ r\geq \zt_0\}$:
\begin{equation}
\ZLA{eqCAnTe:I3TheoremPROPdiR}
\forall\ZEP>0\ \exists \ZSI_\ZEP\,:\ \left (
 \left\{\begin{array}{l}
0<\ZSI<\ZSI_\ZEP\\
t\geq \zt_0\\
r\geq \zt_0
\end{array}\right. \quad \implies |R(t,r;\zt_0-\ZSI)-R(t,r;\zt_0)|<\ZEP
\right )\,.
\end{equation}

Now we prove right continuity at $\zt_0<\zt_M$. We rewrite~(\ref{eqA:I3TheoremPROPdiR})  as the following Fredholm integral equation in the variable $t\in [\zt_0+h,T]$:
\begin{multline*}
\Delta(t,r;h)=\int _{\zt_0}^{\zt_0+h} K(t,\nu)B^* R(\nu,r;\zt_0)\ZD\nu\\
-\int _{\zt_0+h}^ T K(t,\nu)B^* \Delta(\nu,r;h)\ZD\nu
\end{multline*}
so that
\begin{multline}
\ZLA{eqC:I3TheoremPROPdiR}
 \Delta(t,r;h)=\int _{\zt_0}^{\zt_0+h} K(t,\nu)B^* R(\nu,r,\zt_0)\ZD\nu\\
 -\int _{\zt_0+h} ^T R(t,s;\zt_0+h)
  \left [\int _{\zt_0}^{\zt_0+h} K(s,\nu)B^* R(\nu,r;\zt_0)\ZD\nu \right ]
 \ZD s\,.
\end{multline}
The boundedness property in the   Statement~\ref{I1TheoremPROPdiR} shows that
the right hand side tends to zero for $h\to 0^+$ and we have uniform right continuity at $\zt_0\in [0,\zt_M)$, in the following sense  
\begin{equation}
\ZLA{eqCbis:I3TheoremPROPdiR}
\forall\ZEP>0\ \exists h_\ZEP>0\,: \ \left  ( \left\{\begin{array}{l} 
0<h<h_\ZEP\,,\\
\zt_0<t\\
\zt_0<r
\end{array}\right.\   \implies\ 
|R(t,r;\zt_0+h)-R(t,r;\zt_0)|<\ZEP   \right )\,.
\end{equation}

In conclusion, $\zt\mapsto R(t,s;\zt)$ is continuous at every $\zt\in [0,\min\{t,\,r\}]$.
\item[\bf Differenziability:] 
We prove left differentiability. 
We divide with $-\ZSI$ both the sides of~(\ref{eqB:I3TheoremPROPdiR}) and we compute the limit for $\ZSI\to0^+$. The boundedness property in the 
Statement~\ref{I1TheoremPROPdiR},  Lipschitz property~(\ref{eq:LipContRunizt}) and left continuity of $\zt\mapsto R(t,r;\zt)$ uniform respect to $(t,r)$ imply
\[
\lim_{\ZSI\to 0^+}\left [\dfrac{1}{-\ZSI}\int _{\zt_0-\ZSI}^{\zt_0} K(t,\nu)B^* R(\nu,r;\zt_0-\ZSI)\ZD\nu\right ]=-K(t,\zt_0)B^*R(\zt_0,r;\zt_0)
\]
so that ($\partial_{-, \zt}$ denotes the left derivative at $\zt$)
\begin{align*}
\partial_{-, \zt}R(t, r ,\zt_0)=K(t,\zt_0)B^*R(\zt_0,r;\zt_0)\\
-\int_{\zt_0}^T R(t,s;\zt_0)K(s,\zt_0)B^* R(\zt_0,r;\zt_0)\ZD s\,.
\end{align*}
This equality shows that 
\[
\zt_0\mapsto \partial _{-, \zt}R(t, r ,\zt_0)
\]
is continuous on the interval $[0,\zt_M]$. The proof is finished because a continuous function with continuous left (or right) derivative is of class $C^1$.\zdia
\end{description}

Now we examine in more details the function $p_\zt(t)$ in~(\ref{eq:FredhRISOLTA}).
Both the addenda in~(\ref{eq:FredhRISOLTA}) depend linearly and continuously  on $(\hat\xi,\tilde\xi,y)\in M^2_\zt\times L^2(\zt,T)$. So, $p$ is the sum of three terms, which depend separately on $\hat\xi\in\zzr^d$, $\tilde\xi\in L^2(0,\zt;\zzr^d)$ and on the reference signal $y\in L^2(\zt,T;\zzr^p)$. 
We study separately the contribution to the right hand side of~(\ref{eq:FredhRISOLTA})  of $\hat\xi$, $\tilde\xi$ and $y$.
 
\paragraph{\bf The contribution of $\hat\xi$:} the contribution of $\hat\xi $ is
\begin{multline*}
Q_0(t;\zt)\hat\xi=\left [\int_t^TZ (s-t)C^*CZ^*(s-\zt)\ZD s\right.\\
\left. -\int_\zt^T R(t,r;\zt)
\int_r^TZ (s-r)C^*CZ^*(s-\zt)\ZD s\, 
\ZD r
 \right ]\hat\xi \,.
\end{multline*}
The first addendum is of class $C^1$ on the angle $t\geq\zt$ thanks to the regularity properties of $Z(t)$. The second addendum is of class $C^1$ too. Continuity of the derivative respect to $t$  follows by Statement~\ref{I22TheoremPROPdiR} of Theorem~\ref{TheoremPROPdiR} and  Lemma~\ref{lemmaCONTIinteg} applied to
\[
R(t,r;\zt)
\int_r^TZ (s-r)C^*CZ^*(s-\zt)\ZD s\,.
\]

while continuity of the derivative respect to $\zt$ follows from the Statement~\ref{I3TheoremPROPdiR} of Theorem~\ref{TheoremPROPdiR}.

The fact that $\zt\leq t\leq T$ and~(\ref{EqCONDIfinR}) imply $Q_0(T;\zt)=0$ for every $\zt\in[0,t]$. In particular, $Q_0(T,T)=0$.

\paragraph{\bf The contribution of $\tilde\xi$:} We note that
 \begin{multline*}
 \int_t^T Z(s-t)C^*C\int_\zt^s Z^*(s-r)\int_0^\zt N(r-\nu)\tilde\xi(\nu)\ZD\nu\,\ZD r\,\ZD s\\
 =\int_0^\zt \underbrace{\int_t^T Z(s-t)C^*C\int_\zt^s Z^*(s-r)N(r-\nu)\ZD r\,\ZD s
 }_{Q_{1,a}(t,\nu;\zt)}\,\tilde\xi(\nu)\ZD \nu
 \end{multline*}
The function 
\[
Q_{1,a}(t,\nu;\zt)=\int_t^T Z(s-t)C^*C\int _{\zt-\nu}^{s-\nu} Z(s+\nu-\mu)N(\mu)\ZD\mu
\]  
is of class $C^1$ (in the variable $(t,\nu,\zt)\in\{\zt\in[0,T]\,,\ 0\leq t, \nu\leq\zt\}$)   and
\[
\mbox{$Q_{1,a}(T,\nu,\zt)=0$ (in particular $Q_{1,a}(T,\nu;T)=0$)}\,.
\]
Note that if $\nu=T$ it must be $\zt=T$ and so we have also $Q_{2,a}(T,T;T)=0$.

The function $\tilde\xi$ appears also in the addendum
\begin{align*}
-\int_\zt^T R(t,r;\zt)\int_0^\zt Q_{1a}(r,\nu;\zt)\tilde\xi(\nu)\ZD\nu\,\ZD r\\
=-\int_0^\zt\underbrace{ \int_\zt^T R(t,r;\zt)Q_{1a}(r,\nu;\zt)\ZD r  }_{Q_{1b}(t,\nu;\zt)}\, \tilde\xi(\nu)\ZD\nu\,.
\end{align*} 
The functions $Q_{1,b}(t,\nu;\zt)$ is of class $C^1$ in the variable $ (t,\nu,\zt)$. 
This is seen by invoking Statement~\ref{I2TheoremPROPdiR} of Theorem~\ref{TheoremPROPdiR} and Lemma~\ref{lemmaCONTIinteg}.

If either $t$, $\nu$ or $\zt$    equals  $T$
then the function $Q_{1,b} $ is equal zero.

In conclusion, the contribution of $\tilde \xi$ is
\[\int_\zt^T Q_1(t,\nu;\zt)\tilde\xi(\nu)\ZD\nu\,,\qquad
Q_1(t,\nu;\zt)=Q_{1,a}(t,\nu;\zt )-Q_{1,b}(t,\nu;\zt)\ 
\]
and $Q_1(t,\nu;\zt)$ is of class $C^1$  (in the variable $(t,\nu,\zt)$).   
  
 
 \paragraph{\bf The contribution of reference signal $y$:}
 the reference signal appear in the following sum:
 \begin{multline*}
- \int_t^T Z(s-t)C^* y(s)\ZD s+\int_\zt^T \left[\int_\zt^s R(t,r;\zt)Z(s-r)C^*\ZD r\right] y(s)\ZD s\\
=-\int_\zt ^T  \underbrace{ \heaviside(s-t)Z(s-t)C^*}_{=Q_{2,a}(t,s;\zt)}\,y(s)\ZD s+
\int_\zt^T \underbrace{  \int_\zt^s R(t,r;\zt)Z(s-r)C^*\ZD r}_{=Q_{2,b}(t,s;\zt)}\, y(s)\ZD s\,.
 \end{multline*}
We note an abuse of notations in the definition of $Q_{2,a}$: the function $Z(r)$ \emph{is not defined if $r<0$.  We wrote $\heaviside(r) Z(r)$ to intend that  we put \emph{by definition} $\heaviside(r)Z(r)=0$  when $r<0$.} So,
$r\mapsto\heaviside(r)Z(r)$ has a jump at $r=0$. 

The contribution of $y$ is
\[
\int _\zt ^T Q_{2} (t,s;\zt) y(s)\ZD s\,,\qquad Q_2(t,s;\zt)=-Q_{2a}(t,s;\zt)+Q_{2b}(t,s;\zt)\,.
\]

We can state  the following properties. Let
\begin{equation}\ZLA{eq:DEFIaTAU}
A _{\zt}= \{0\leq \zt<T\,,\ \zt\leq t\leq T\,,\ \zt\leq s\leq T\}\,.
\end{equation}
\begin{enumerate}
\item
The functions $Q_2(t,s,\zt)$, $\partial_t Q_2(t,s,\zt)$, $\partial_s Q_2(t,s,\zt)$ are continuous as functions of $(t,s,\zt)$ on the set $A_\zt\setminus \{(y,y)\}$ with finite limits for $(t,s,\zt)\to (y,y,\zt)$.
\item for $t$ and $s$ fixed (even when $s=t)$ the function $\zt\mapsto Q_2(s,t,\zt)$ is continuous, even differentiable, and
  the function $(s,t,\zt)\mapsto \partial_\zt Q_2(s,t,\zt)$ is continuous on $A_\zt$.
  \item if $\zt=T$ then $Q_2(s,t,T)=0$.
\end{enumerate}

 \paragraph{The form of $\upzt$ and of $\wpzt$}

 We recall that the optimal control is $\upzt=-B^*p(t)$ so that 
the transformation
\[
(\hat\xi,\tilde\xi,y)\mapsto \upzt
\] 
 is linear and continuous from 
  $(\hat\xi,\tilde\xi,y)\in M^2_\zt\times C([\zt,T];\zzr^d)$  to 
$C([\zt,T],\zzr^d)$, explicitly given by
 \begin{subequations}
\begin{equation}\ZLA{eqOptiContdopoFred}
\upzt(t;\Xi_\zt)=-\left [B^*Q_0(t;\zt)\hat\xi+
\int_0^\zt B^*Q_1(t,s;\zt)\tilde\xi(s)\ZD s+\int_\zt^T B^*Q_2(t,s;\zt) y(s)\ZD s
\right ]\,.
\end{equation}
  \begin{Remark}\ZLA{RemaSULLAformSEMIfeedoptimContDA FREDH}{\rm
 The previous formulas hold for $t\geq \zt$. In particular, by taking the limit for $t\downarrow \zt$, we see that $\upzt(\zt;\Xi_\zt)  $ is a linear functional of $(\hat\xi,\tilde\xi\ZCD)\in M^2_\zt$ and   of the restriction of the reference signal $y$ to the future interval of time $[\zt,T]$.

A sharper form of this observation depends on the next Lemma~\ref{lemmaSemiEVOLpropeOPTIMAL}   and it is given in Remark~\ref{RemaSULLAformSEMIfeedoptimContDA FREDHbis}.\zdia 
 }\end{Remark}

We replace in~(\ref{eq:FormVARconst}) and we find
\begin{equation}\ZLA{eqOptiwWwdopoFred}
\wpzt(t;\Xi_\zt)=H_0(t;\zt)\hat\xi+\int_0^\zt H_1(t,s;\zt)\tilde\xi(s)\ZD s+\int_\zt^T H_2(t,s;\zt)y(s)\ZD s
\end{equation}
\end{subequations}
where:

\begin{align*}
H_0(t;\zt)&=Z^*(t-\zt)-\int_\zt^t Z^*(t-r)BB^*Q_0(r;\zt)\ZD r\,,\\
H_1(t,s;\zt)&=\int_\zt^t Z^*(t-r)N(r-s)\ZD r-\int_\zt^t Z^*(t-r)BB^* Q_1(r,s;\zt)\ZD r\,,\\
H_2(t,s;\zt)
&=-\int_\zt ^t Z^*(t-r)BB^*Q_2(r,s;\zt)\ZD r\,.
\end{align*}
The matrices $H_i$ have the following properties.
\begin{description}
\item[\bf The matrix $H_0(t;\zt)$] is of class $C^1$ in the angle $t\geq \zt$ and
$
H_0(T,T)=I
$.
 
\item[\bf The matrix $H_1(t,s;\zt)$]  is of class $C^1$ in the variable $(t,s,\zt)$. The regularity of the second integral follows from the regularity of $Q_1(r,s;\zt)$. The regularity of the first integral is seen by substituting $r=\nu+s$.

If $\zt=T$ then $H_1=0$.
\item[\bf The matrix $H_2(t,s;\zt)$.] Let $A_\zt$ be the set in~(\ref{eq:DEFIaTAU}). The functions $H_2(t,s,\zt)$, $\partial_s  H_2(t,s,\zt)$,
$\partial_\zt  H_2(t,s,\zt)$ are continuous on $A_\zt$ while the function $\partial_t  H_2(t,s,\zt)$ is continuous on $A_\zt$, provided that $s\neq t$.

If $\zt=T$ then $H_2(t,s;T)=0$.

 \end{description}
 
 We close this section by mentioning that the use of Fredholm integral equations for systems with \emph{finite} delays was  advocated in~\cite{ManitiusOPTIMALtrieste76}.
\section{\ZLA{sectTheDI}The dissipation inequality and its consequences}

In this section we use Eq.~(\ref{eq:SisteVoltezt}) but with different initial times, 
say $\zt$ and $\ztu$. So, by using  the notation~(\ref{eq:FormVARconst}), the
 solution is  respectively
  
\[
r\mapsto w (r ;\zt, \Xi_\zt,u)\,,\quad\ r\mapsto w (r;\ztu,\Xi_{\ztu}, u)\,.
\]
We intend  respectively $r\in [\zt,T]$ or $r\in[\ztu,T]$,
   $\Xi_\zt\in M^2_\zt$, or 
$\Xi_{\ztu} \in M^2_\ztu$ and we intend  that $u$ is defined in the proper interval, in the interval  $[\zt,T]$ in the first case and in the interval $[\ztu,T]$ in the second case.

For $r\in [\zt,T]$ we define
\[
\Xi (r;\zt,\Xi_\zt,u )=\left (w (r;\zt, \Xi_\zt,u ),  
  w ( s;\zt, \Xi_\zt,u  )_{|_{  (0,r)}}  
     \right ).
\]
Keep in mind~(\ref{eq:DefiESTEdiW})   in order to interpret the restriction of $w (\cdot;\zt, \Xi_\zt,u  )$ to $(0,r)$.

 We retain the notations $\upzt(\cdot;\Xi_\zt)$ and $\wpzt(\cdot;\Xi_\zt)$ for the optimal control and the corresponding solution of~(\ref{eq:SisteVoltezt}):
 \[
 \wpzt(\cdot;\Xi_\zt)=w(\cdot;\zt,\Xi_\zt,\upzt(\cdot;\Xi_\zt))\,.
 \]
 
 We introduce the notation
 \[
\Xi^+_\zt( r;\Xi_\zt)=
\left (w^ +_\zt(r; \Xi_\zt  ),
  (w^+ _\zt(\cdot; \Xi_\zt  ) _{|_{[0,r]}}      
     \right ). 
\]
 

We recall that the optimal control is continuous so that   
\begin{multline*}
W_\zt(\Xi_\zt )=_{\rm\small def} \min_{u\in L^2(\zt,T;\zzr^m)}J_\zt(\Xi_\zt ;u)
\\=\int_\zt^T\left [ \|C \wpzt(r;\Xi_\zt)-y(r)\|^2+\|\upzt(r;\Xi_\zt)\|\right ]\ZD r=
 \min_{u\in C([\zt,T];\zzr^m)}J_\zt(\Xi_\zt;u)\,.
\end{multline*}
So, when studying the properties of the function $W_\zt$ we can confine ourselves to use continuous or piecewise continuous controls.
 
We fix an initial condition $\Xi_\zt$ at the initial time $\zt$.  
Let $\ztu\in(\zt,T]$.   Let $u \in C([\zt,\ztu];\zzr^m)$. We consider the solution $r\mapsto w (r;\zt,\Xi_\zt,u)$ of~(\ref{eq:SisteVoltezt}) on $[\zt,\ztu]$.

We extend the control $u $ to a control $u_1$ on $[\ztu,T]$ and we solve the   equation~(\ref{eq:SisteVoltezt}) on the entire interval $[\zt,T]$, with the same initial condition $\Xi_\zt$ and  with the control
\[
u_{e}(r)=\left\{\begin{array}{lll}
u (r)&{\rm if}& r\in [\zt,\ztu)\\
u_1(r)&{\rm if}& r\in [\ztu,T)\,.
\end{array}\right.
\]

It is clear that
\[\left\{\begin{array}{l}
w (r; \zt, \Xi_\zt, u_{e}) =w (r; \zt, \Xi_\zt, u )\qquad r \in [\zt ,\ztu]\\
w  (r;\zt,  \Xi_\zt, u_{e}) =w(r;{\ztu}, \Xi (\ztu;\zt, \Xi_\zt,u ), u_1)\qquad r\in[\ztu ,T]\,.
\end{array}\right.
\]
We use Bellman optimality principle: we have
\begin{multline}\ZLA{eq:PREperCONCAcontOtti}
W_\zt(\Xi_\zt)\leq
J_\zt(\Xi_\zt,u_{e})= \int_\zt^{\ztu}  \left [
 \|C  w  (r;\zt, \Xi_\zt,u )-y(r)\|^2+\|u (r)\|^2
\right ] \ZD r  \\   +J_{\ztu}( \Xi (\ztu;\zt,\Xi_\zt,u);u_1)\,.
\end{multline}
By taking the minimum of the right hand side respect to $u_1$, we get
\begin{equation}\ZLA{DIformaintegr}
W_\zt(\Xi_\zt)\leq \int_\zt^{\ztu}  \left [
 \|C  w  (r;\zt, \Xi_\zt,u )-y(r)\|^2+\|u (r)\|^2 
\right ] \ZD t + W_{\ztu} ( \Xi (\ztu;\zt,\Xi_\zt,u ))\,.
\end{equation}
This is the form taken by the \emph{dissipation inequality (DI) in integral form} in the presence of the reference signal $y$. 

The definition of $W_\zt$ and the unicity of the optimal control shows that\emph{~(\ref{DIformaintegr}) holds as an equality if and only if}
\[
u(r)=u^+_\zt(r;\Xi_\zt) \quad{\rm if}\ r\in [\zt,\ztu)\,.
\]
and
 \begin{Lemma}\ZLA{lemmaSemiEVOLpropeOPTIMAL}
When the initial condition at the initial time $\ztu$ is $\Xi^+_\zt(\ztu;\Xi_\zt)$ then the optimal 
control on $[\ztu,T]$ is $\left (\upzt(\cdot;\Xi_\zt)\right )_{|_{[\ztu,T]}}$:
 \begin{align*}
   u^+_{\ztu} (r;\Xi^+_{\zt }( \ztu;\Xi_\zt)) =\upzt(r,\Xi_\zt)\,,  &\quad   r\in([\ztu,T]\,,\\
   w^+_{\ztu} (r;\Xi^+_{\zt }( \ztu;\Xi_\zt)) =\wpzt(r,\Xi_\zt)\,,  &\quad    r\in([\ztu,T] \,.
   \end{align*}
 
 \end{Lemma}

We divide both the sides of~(\ref{DIformaintegr}) with $\ztu-\zt$ and we rewrite it as
\begin{multline*}
-
\dfrac{1}{\ztu-\zt}\int_\zt^{\ztu}  \left [
 \|C  w (r;\zt,\Xi_\zt,u )-y(r)\|^2+\|u (r)\||^2 
\right ] \ZD t\\ \leq \dfrac{ W_{\ztu} ( \Xi  ({\ztu};\zt,\Xi_\zt,u ))-W_\zt(\Xi_\zt)}{\ztu-\zt}\,.
\end{multline*}
This inequality holds for every continuous control $u$.
We compute   the limit of both the sides for   $\ztu\downarrow\zt$. We get the \emph{dissipation inequality (DI) in differential form}
 \begin{equation}
 \ZLA{DIforaDiffe}
 -\left [ \|C \hat\xi-y(\zt)\|^2+\|u(\zt)\|^2\right ]\leq \left ( D_+ W_{\ztu} ( \Xi (\ztu;\zt,\Xi_\zt,u ))\right )_{|_{t=\zt}} 
 \end{equation}
 where $D_+$ denote the first right Dini derivative respect to $t$:
 \[
D_+ f(t)=\liminf _{h\to 0_+} \dfrac{f(t+h)-f(t)}{h} \,.
 \]
 
   \emph{We shall see that  $W_{\ztu} ( \Xi (\ztu;\zt,\Xi_\zt,u ))\in C^1([\zt,T])$ so that the Dini derivative in~(\ref{DIforaDiffe})  is just the usual right derivative and, in particular,    }
 \begin{equation}
 \ZLA{eq:DiniDeriComeLimiDeri}
 \left ( D_+ W_{\ztu} ( \Xi (\ztu;\zt,\Xi_\zt,u ))\right )_{|_{\ztu=\zt}}  
=\lim _{\ztu\to\zt^+}\dfrac{\ZD}{\ZD \ztu}\left (   W_{\ztu} ( \Xi (\ztu;\zt,\Xi_\zt,u ))\right )\,.
 \end{equation}
 
 The inequality~(\ref{DIforaDiffe}) holds as an equality if $u(t)=u^+(t)$ and so we can state:  
\begin{equation}
\ZLA{eq:PreClosOptiContro}
\upzt(\zt;\Xi_\zt)=\arg\min\left \{
|u(\zt)|^2+
 \left ( D_+ W_{\ztu} ( \Xi (\ztu;\zt,\Xi_\zt,u ))\right )_{|_{t=\zt}}
\right \}\,.
\end{equation}

 The derivative in~(\ref{eq:PreClosOptiContro}) is the derivative of the composite function: that of $W_\ztu$ computed along $\Xi(\cdot;\zt,\Xi_\zt, u)$, driven by the continuous control $u$.
 
By explicitly computing the derivative in~(\ref{eq:DiniDeriComeLimiDeri}), we shall see that the minimization problem in~(\ref{eq:PreClosOptiContro}) admits a unique solution so that $u=u^+(\zt;\Xi_\zt)$ if and only if it realizes the equality in~(\ref{DIforaDiffe}).  After that we can use~(\ref{eq:PreClosOptiContro}) to derive an expression of the optimal control which is more convenient then~(\ref{eqOptiContdopoFred}). This expression will depend on certain coefficients and we must find equations for them. We are going to find a set of differential equations for the coefficients which in particular includes a version of the Riccati equation  for systems with memory. This program is realized in the next section.

 \subsection{\ZLA{sec:FinalRepre}The sinthesis of the optimal control}
 
Our goal in this section is the solution of the optimization problem~(\ref{eq:PreClosOptiContro}) and the derivation of a new expression of the optimal control, which do not require the solution of a Fredholm integral equation. 
This  goal is achieved in the next steps.

\paragraph{Step~1: the function  $  W_{\ztu} (\cdot )$.}
We fix any $\ztu\in [\zt,T]$ and an \emph{arbitrary} initial condition in $M^2_\ztu$. In the next steps the  initial condition will be specialized and will be taken equal to $ \Xi (\ztu;\zt,\Xi_\zt,u )$. But, in this step it is an \emph{arbitrary} element of $M^2_\ztu$. So, in order to avoid confusions, the \emph{ arbitrary } initial condition in $M_\ztu^2$ which we use now is denoted $\ZOMq_\ztu$:
\[
\ZOMq_\ztu=(\hat\ZOM,\tilde\ZOM\ZCD)\in M^2_\ztu\,.
\]
 
 We must compute 
\[
W_{\ztu}(\ZOMq_\ztu)= \int _\ztu ^T \left  [ \| Cw^+_\ztu(r; \ZOMq_\ztu  ) -y(r)\|^2+\|u^+_\ztu(r; \ZOMq_\ztu  )\|^2\right  ]\ZD r
\]
The functions $w^+_\ztu(r; \ZOMq_\ztu  )$ and $u^+_\ztu(r; \ZOMq_\ztu  )$
are given by~(\ref{eqOptiContdopoFred}) and~(\ref{eqOptiwWwdopoFred})
(with the obvious change of the notations, $t$ replaced by $r$, $\zt$   by $\ztu$ and $\xi$ by $\ZOM$).
 
We replace~(\ref{eqOptiContdopoFred}) and~(\ref{eqOptiwWwdopoFred}). 
We have\footnote{note that  we are working in \emph{real spaces} but  we use the notation
\[
2\zreal \ZL a,b\ZR \quad \mbox{to denote}\quad  \ZL a, b\ZR+\ZL b, a\ZR\,.
\]}:
 
\begin{multline*}
 \| Cw^+_\ztu(r;\ZOMq_\ztu)-y(r)\|^2=\ZL \hat  \ZOM,H_0^*(r;\ztu   )C^*CH_0(r; \ztu  )\hat  \ZOM\ZR\\
 +2\zreal\left   \{
 \ZL\hat  \ZOM,H_0^*(r;\ztu  )C^*C\int  _0^\ztu   H_1(r,s;\ztu  )\tilde \ZOM(s)\ZD s\ZR\right.\\ 
\hskip 2cm \left. +
  \ZL\hat  \ZOM,H_0^*(r;\ztu  )C^*C\int  _\ztu^T   H_2(r,s;\ztu  )y(s)\ZD s\ZR
 \right   \}\\
\hskip 2cm -2\zreal\ZL \hat  \ZOM,H_0^*(r;\ztu  )C^* y(r)\ZR\\
 +\int  _0^\ztu  \int  _0^\ztu  \ZL\tilde \ZOM(s),H_1^*(r,s;\ztu  )C^*CH_1(r,\nu;\ztu  )\tilde \ZOM(\nu)\ZR \ZD \nu\,\ZD s\\
 +2\zreal \left   \{\int  _0^\ztu  \int  _\ztu  ^T \ZL\tilde \ZOM(s), H_1^*(r,s;\ztu )C^*CH_2(r,\nu;\ztu  )y(\nu)\ZR\ZD\nu\ZD s\right.\\
   \hskip .5cm    \left. 
 - \int  _0^\ztu   \ZL  \tilde \ZOM(s),H_1^*(r,s;\ztu  ) C^*y(r)  \ZR  \ZD s 
 \right   \} +M(r,\ztu )
   \\
\end{multline*} 
where 
\[ 
M(r,\ztu ) =\left   \|\int  _\ztu^T   CH_2(r,s;\ztu  )y(s)\ZD s\right   \|^2
+\|y(r)\|^2  
-2\zreal   \int  _\ztu  ^T\ZL CH_2(r,s;\ztu  )y(s),y(r)\ZR\ZD s  \,.
\]   
Note that  $M(r,\ztu ) =0$ when $y=0$ and that
\begin{itemize}
\item $H_0^*(r;\ztu   )C^*CH_0(r; \ztu  )$ is selfadjoint;
\item $\left [H_1^*(r,s;\ztu  )C^*CH_1(r,\nu;\ztu  )\right ]^*=
H_1^*(r,\nu;\ztu  )C^*CH_1(r,s;\ztu  )$. So, the computation of the adjoint exchanges the roles of $s$ and $\nu$.
\end{itemize}

The expression of $\|u(r)\|^2$ is similar, with the matrices $CH_i$ replaced by the corresponding matrices $-B^*Q_i$ and $y=0$. 
We sum the two contributions and we integrate in the variable $r$ on the 
interval  $[\ztu ,T]$ (we recall $\zt  <\ztu <T$. We find
\begin{multline}
\ZLA{eq:La formaFINlediW}
 W_{\ztu  }(\ZOMq_{\ztu } )=\ZL \hat\ZOM,P_0( \ztu  )  \hat\ZOM\ZR+2\zreal \ZL \hat\ZOM,\int  _0^\ztu   P_1(s,\ztu  )\tilde\ZOM(s)\ZD s\ZR\\
+\int  _0^\ztu  \int  _0^\ztu   \ZL \tilde\ZOM (s),P_2(s,\nu,\ztu  ) \tilde\ZOM(\nu)\ZR\ZD\nu\,\ZD s\\
+2\zreal\ZL \hat\ZOM,\dd(\ztu )\ZR+2\zreal \int  _0^\ztu   \ZL \tilde\ZOM(s),\D(s,\ztu)\ZR\ZD s 
+M(\ztu  )\,.
\end{multline}
Explicit expressions of the matrices $P_i$, of the vectors $d_i $ and of $M$   in~(\ref{eq:La formaFINlediW}) are easily derived. We confine ourselves to state the properties which we are going to use.
 \begin{enumerate}
 \item the matrices $P_i$ depend neither on the initial condition $(\hat\ZOM,\tilde\ZOM\ZCD)$ nor on the reference signal $y$ and:
 \begin{enumerate}
 \item the following equalities hold:
\[
P_0(\ztu)=P_0^*(\ztu)\,,\qquad P_2(s,\nu,\ztu)=P_2^*(\nu,s,\ztu)\,.
\]
\item $P_0\ZCD\in C^1([0,T])$ and $P_0(T)=0$.
\item the function $P_1(s,t)$ is defined on the set $\{(s,t)\,,\ 0\leq s\leq t\leq T\}$. On this set both $P_1(s,t)$ and $\partial_t P_1(s,t)$ are continuous and $P_1(s,T)=0$ for every~$s$.
\item $P_2(s,\nu,t)$ is defined on the set $\{(s,\nu,t)\,,\ 0\leq s\,,\, \nu\leq t\leq T\}$. On this set, $P_2(s,\nu,t)$ and $\partial_t P_2(s,\nu,t)$ are continuous. Moreover, $P_2(s,\nu,T)= 0$ for every $s$ and $\nu$.
\end{enumerate}
\item the vectors $d_i$ and $M$ do not depend on the initial condition $(\hat\ZOM,\tilde\ZOM\ZCD)$. They depend on the reference signal $y$. If $y=0$ then we have $d_i=0$ and $M=0$. Furthermore:
\begin{enumerate}

\item  $\dd(t)$ and $M(t)$ are of class $C^1$ on the interval $[0,T]$ and $\dd(T)=0$, $M(T)=0$.
\item  $\D(s,t)$ is defined on the triangle $\{0\leq s\leq t\leq T\}$. On this triangle both $  \D(s,t)$ and $ \partial_t \D(s,t)$ are continuous. And: $\D(s,T)=0$ for every $s$.  
\end{enumerate}
 
\end{enumerate}

\paragraph{Step~2: the proof of~(\ref{eq:DiniDeriComeLimiDeri}) and the   optimal control.}

We replace $\ZOMq_t$ with $\Xi(t;\zt,\Xi_\zt,u) $  in~(\ref{eq:La formaFINlediW}). As we noted, we can use a continuous control  $u$ so that 
 the composite function so obtained is of class $C^1$: $W_{\ztu} ( \Xi (\ztu;\zt,\Xi_\zt,u ))\in C^1([\zt,T])$ and the equality~(\ref{eq:DiniDeriComeLimiDeri}) holds.

We compute the derivative and its limit for $\ztu\downarrow \zt$; we replace it in~(\ref{eq:PreClosOptiContro}).  Then we compute the minimum respect to $u$. The computation is boring but elementary. So, we confine ourselves to state the results. 

The functional which has to be minimized in~(\ref{eq:PreClosOptiContro})   in order to obtain $u^+_\zt(\zt;\Xi_{\zt})$ is the sum of a quadratic form $\mathcal{F}(\Xi_\zt,y_{[\zt,T]})$,   which does not depend on $u $, plus the the following  quadratic coercive form of $u\in\zzr^m$:
\begin{multline}\ZLA{eq:functDAminimizzare}
\left \|
u+B^*\left [P_0(\zt)\tilde\xi+\int_0^\zt P_1(s,\zt)\tilde\xi(s)\ZD s+d_1(\zt)\right ]\,\right \|^2\\
-\left \|B^*\left [P_0(\zt)\tilde\xi+\int_0^\zt P_1(s,\zt)\tilde\xi(s)\ZD s+d_1(\zt)  \right ]\,\right \|^2
\end{multline}
and so, as we announced, the minimum is reached at a unique vale of $u$, which is value  the optimal control:
\begin{equation}
\ZLA{eq:AlteFORMoptimCONTR}
u^+_\zt(\zt;\Xi_\zt)=-B^*\left [P_0(\zt)\hat\xi+\int_0^\zt P_1(s,\zt) \tilde(s)\ZD s+\dd(\zt)\right ]\,.
\end{equation}

Now we can give a sharper form of the observation in Remark~\ref{RemaSULLAformSEMIfeedoptimContDA FREDH}.

\begin{Remark}\ZLA{RemaSULLAformSEMIfeedoptimContDA FREDHbis}{\rm
Equality~(\ref{eq:AlteFORMoptimCONTR}) holds in particular when the initial time is $\zt_0<\zt$ and $\Xi=(\hat\xi,\tilde\xi)=\Xi^+_{\zt_0}(\zt;\Xi_{\zt_0})$. We consider the special case of the problem~(\ref{eq:SisteVolte0})-(\ref{eq:SisteCOSTO0}), i.e. when the initial time is $ 0$. We recall that in this case $\Xi_0=\hat\xi=w(0)$. We use Lemma~\ref{lemmaSemiEVOLpropeOPTIMAL} and we see that   the optimal control $u^+_0(\zt;w(0))$ at any   time $\zt\geq 0$ as given by~(\ref{eq:AlteFORMoptimCONTR})
is
\begin{equation}
\ZLA{eq:AlteFORMoptimCONTRda0}
u^+_\zt(\zt;\Xi_\zt)=-B^*\left [P_0(\zt)w(\zt)+\int_0^\zt P_1(s,\zt) w(s)\ZD s+\dd(\zt)\right ]\,.
\end{equation}
In this formula $w\ZCD$ is the solution of~(\ref{eq:SisteVolte0}) with the control $u^+_\zt(\cdot;\Xi_\zt)$, and so it is $w^+_0(\cdot;\Xi_0)$: the optimal control is
the sum of a linear functional of $y_{|_{[\zt,T]}}$, i.e. of the \emph{future} of the reference signal plus  a linear   a feedback of the  \emph{state $(w^+_0(\zt;w(0)),w^+_0(\cdot;w(0))_{|_{(0,\zt)}})$ at time $\zt$}.\zdia}
 \end{Remark}

Observe that  we have a representation of the optimal control and of $W^+$ in terms of the same coefficient matrices and vectors $P_i$, $d_i$  and $M$. In order that this representation have its advantages, in the next step we present   equations for the coefficients which do not require the solution of a Fredholm integral equation.

\emph{The representation~(\ref{eq:AlteFORMoptimCONTRda0}) of the optimal control and the equations for the coefficients derives in the next Step~3 constitute the synthesis of the optimal control.}

\paragraph{Step~3: the equations of the coefficients.}
 
We use the following notation: a superimposed dot denotes  derivative\footnote{actually, the derivatives are right derivatives but we recall that a continuous function on an interval whose right derivative is continuous too it is of class $C^1$.} respect to $\zt$:
\[
\dot P_0(\zt)=\dfrac{\ZD}{\ZD \zt} P_0(\zt)\,,\quad \dot P_1(s,\zt)=\dfrac{\partial}{\partial \zt} P_1(s,\zt)\quad {\rm etc.}
\]
 We 
fix any initial time $\zt\in [0,T]$ and any initial condition $\Xi_\zt=(\hat\xi,\tilde\xi\ZCD)\in M^2_\zt $. We replace 
  the optimal control~(\ref{eq:AlteFORMoptimCONTR}) in the DI in differential form~(\ref{DIforaDiffe}), which is identically zero for $\ztu\geq \zt$.
An explicit computation gives:
\begin{multline}\ZLA{eq:laDIeQUALzero}
\ZL \hat\xi,\left [\dot P_0(\zt)+C^*C-P_0(\zt)BB^*P_0(\zt)+A^*P_0(\zt)+P_0(\zt )A\right.\\
 \left.  +P_1(\zt,\zt)+P_1^*(\zt,\zt)\Big]\hat\xi\ZR   \right.\\
+2\zreal 
\ZL\hat\xi,\int_0^\zt\left [\dot P_1(s,\zt)+ P_0(\zt) N(\zt-s)+A^*P_1(s,\zt)\right.\\
\hskip 3cm \left.+P_2(\zt,s,\zt)-P_0(\zt) BB^*P_1(s,\zt)\Big]\tilde\xi(s)\ZR\ZD s   \right. \\
+\left\{\int_0^\zt\int_0^\zt \ZL \tilde\xi(s),\left [
\dot P_2(s,\nu,\zt)+N^*(\zt-s)P_1(\nu,\zt)\right.\right.\\
 \left. \qquad\qquad \qquad \qquad+P_1^*(s,\zt)N(\zt-\nu)-P_1^*(s,\zt)BB^*P_1(\nu,\zt)
\Big]\tilde\xi(\nu)\ZR\ZD\nu\,\ZD s\biggr\}\right.\\
 +2\zreal 
 \ZL\hat\xi,\left [\dot \dd(\zt)+A^*\dd(\zt)+\D(\zt,\zt)-P_0(\zt)BB^*\dd(\zt)-C^*y(\zt)\right ]\ZR
\\
+\int_0^\zt \ZL\tilde\xi(s),\left [\dot \D(s,\zt)+N^*(\zt-s)\dd(\zt)-P_1^*(s,\zt)BB^* \dd(\zt)\right ]\ZR \ZD s\\
+\dot M(\zt)+\|y(\zt)\|^2-\|B^*\dd(\zt)\| 
=0\,.
\end{multline}

This equality holds for any arbitrary initial conditions $\Xi_\zt\in M^2_\zt$ and any arbitrary reference signal $y\in L^2(0,T;\zzr^p)$.

We note that the coefficient matrices $P_i$ depend neither on $\Xi $ nor on $y$ while $\dd$, $\D$ and $M$ do depend on $y$ but do not depend on $\Xi$.

Furthermore, $\Xi=(\hat\xi,\tilde \xi)$ and $\hat\xi$ and $\tilde \xi$ are independent. 

We consider the equality~(\ref{eq:laDIeQUALzero}) when $y=0$. The fact that $\hat\xi$ and $\tilde\xi $ can be independently assigned imply that the coefficient matrices $P_i$ solve the following system of equations
\begin{subequations}
\begin{equation}
\ZLA{Eq:EquqRiccati}
\left\{
\begin{array}{l}
\dot P_0(\zt)+A^*P_0(\zt)+P_0(\zt)A+P_1(\zt,\zt)+P_1^*(\zt,\zt)-P_0(\zt)BB^*P_0(\zt)+C^*C=0\\
\dot P_1(s,\zt)+A^*P_1(s,\zt)+P_0(\zt)N(\zt-s)+P_2(\zt,s,\zt)-P_0(\zt)BB^*P_1(s,\zt)=0\\
\dot P_2(s,\nu,\zt)+N^*(\zt-s)P_1(\nu,\zt)+P_1^*(s,\zt)N(\zt-\nu)-P_1^*(s,\zt)BB^*P_1(\nu,\zt)=0\\
\mbox{$P_0(  \zt  )$ defined on $[0,T]$ and $P_0(T)$=0}\\
\mbox{$P_1(s,  \zt  )$ defined on $\{(s,  \zt  )\,,\ 0\leq s\leq   \zt  \leq T$ and $P_1(s,T)=0$}\\
\mbox{$P_2(s,\nu,  \zt  )$ defined on $\{(s,\nu,  \zt  )\,,\ 0\leq s\,,\, \nu\leq   \zt  \leq T\}$ and $P_2(s,\nu,T)=0$}\,. 
\end{array}\right.
\end{equation}
 
The system~(\ref{Eq:EquqRiccati}) is the version of the Riccati equations for systems with memory, first derived in~\cite{PandolfiVOLTERRAieee18} in the special case $A=0$ and $y=0$.

By considering the case $\Xi=0$ and $y\neq 0$ we derive the following equation for $M(t)$:
\begin{equation}
\ZLA{eq:SISTEfinPERcoeffM}
\dot M(\zt)=\|B^* \dd(\zt)\|^2-\|y(\zt)\|^2\,,\qquad M(T)=0\,.
\end{equation}

Finally, by using the fact that    $\Xi$  and $y$ are arbitrary we get the system of  equations
\begin{equation}
\ZLA{LeEquaSINTESdD}
\left\{
\begin{array}{ll}
\dot \dd(\zt)+(A^*- P_0(\zt)BB^*)\dd(\zt)+\D(\zt,\zt)-C^*y(\zt)=0\,,&\  \dd(T)=0\\
\dot\D(s,\zt)+N^*( \zt-s)\dd(\zt)-P_1^*(s,\zt)BB^*\dd(\zt)=0\,,& \D(s,T)=0\,.
\end{array}\right.
\end{equation}
The function $\dd(  \zt  )$ is defined on $[0,T]$ while the function $\D(s,  \zt  )$ is defined on $\{(s,  \zt  )\,,\ 0\leq s\leq   \zt  \leq T\}$.
 
\end{subequations}

These equations extend to systems with persistent memory the sysnthesis of the optimal control.

We end this section by noting that   the synthesis of the optimal control can be described by equations with different forms of those in~(\ref{Eq:EquqRiccati})-(\ref{LeEquaSINTESdD}). In fact,~(\ref{eq:La formaFINlediW}) can be written in several equivalent manners, as for example
\begin{multline*}
W_{\zt   }(\ZOMq_{\zt  } )=\ZL \hat\ZOM,P_0( \zt   )  \hat\ZOM\ZR+2\zreal \ZL \hat\ZOM,\int  _0^\zt    P_1(s,\zt -s  )\tilde\ZOM(\zt -s)\ZD s\ZR\\
+\int  _0^\zt   \int  _0^\zt    \ZL \tilde\ZOM (\zt -s),P_2(\zt-s,\zt-\nu,\zt   ) \tilde\ZOM(\zt - \nu)\ZR\ZD\nu\,\ZD s\\
+2\zreal\ZL \hat\ZOM,\dd(\zt  )\ZR+2\zreal \int  _0^\ztu   \ZL \tilde\ZOM(s),\D(\zt -s,\zt )\ZR\ZD s 
+M(\zt   )
\end{multline*}
or, we can rename the matrices $P_i$ in~(\ref{eq:La formaFINlediW}) as follows
\[
\tilde P_0(\zt )=P_0(\zt )\,,\quad \tilde P_1(s,\zt )=P_1(\zt -s,\zt )\,,\quad \tilde P_2(s,\nu,\zt)=P_2(\zt-s,\zt-\nu,\zt)
\]
 (and, possibly, $\tilde \D(s,\zt)=\D(\zt-s,\zt)$).
 
 Computations similar to the ones above can be done, leading to different forms of the equations~(\ref{Eq:EquqRiccati})-(\ref{LeEquaSINTESdD}). We leave the computations to the interested readers.
 
 \section{\ZLA{OpeFORMequa}A state space representation of the system and of the synthesis of the optimal control}
 
 We represent the equations~(\ref{Eq:EquqRiccati}) and~(\ref{LeEquaSINTESdD}) in operator form, by introducing the following operators  which are then used to represent the evolution   of the state $t\mapsto \Xi(t;0,\Xi_0,u)$ of the system with memory. We represent the memory term as
 \begin{equation}\ZLA{secSTATEequaConvMEM}
\intt N(t-s)w(s)\ZD s=\intt N(s)\tilde \xi(s,t)\ZD\,s\qquad \tilde\xi(s)=w(t-s)\,. 
 \end{equation}
 This way, the pair $t\mapsto\Xi(t;0,\Xi_0,u)=(w(t),\tilde\xi(t))$ solves
 \begin{equation}\ZLA{eq:secSTATEstarerepreS}
 \begin{array}{ll}
\displaystyle  D_t w(t)=Aw(t)+\intt N(s)\tilde\xi(s,t)\ZD s+Bu(t)\,,\quad w(0)=\hat\xi_0 & \mbox{(for $t>0$)}\\[2mm]
\displaystyle 
D_t \tilde\xi(s,t)=-D_s \tilde \xi(s,t)\,,\quad  \tilde\xi(0,t)=w(t)&
\mbox{(for $0<s<t$)}
\end{array} 
 \end{equation}
 (we use $D_t $ and $D_s$ to denote the derivative respect to the indicated variables). 
 
 Let $\zt$ be any time, $\zt\in (0,T)$.
 We denote $\ZNN_\zt$ the memory operator in~(\ref{secSTATEequaConvMEM}), written in terms of $\tilde\xi$:
\[
\ZNN_\zt\in\mathcal{L}(L^2(0,\zt;\zzr^d),\zzr^d)\,,\quad \ZNN_\zt \tilde\xi =\intzt N( s)\tilde\xi(s)\ZD s
\] 
and,  consistent with this representation of the memory   we introduce
 $\PP_1(\zt)$: $L^2(0,\zt)\mapsto \zzr^d$ and $\PP_2(\zt)$: $L^2(0,\zt)\mapsto L^2(0,\zt)$:
\[
\PP_1(\zt)\tilde\xi=\intzt P_1(\zt-s,\zt)\tilde \xi(s)\ZD s\,,\qquad\left ( \PP_2(\zt)\tilde\xi\right )(s)=\intzt P_2(\zt-s,\zt-\nu,\zt)\tilde \xi(\nu)\ZD\nu\,.
\]
Note that $ \PP_1^*(\zt)\,:\ \zzr^d\mapsto L^2(0,\zt) $ is the operator
 \[
 \left (\PP_1^*(\zt)\hat\xi\right )\ZCD=P_1^*(\zt-\cdot,\zt)\hat\xi\,.
 \]
 The \emph{Riccati operator} is
  \[
\PP ( \zt)=\left[
\begin{array}{cc}
P_0(\zt)& \PP_1 (  \zt)\\
\PP_1^*(  \zt)&\PP _2( \zt)
\end{array}
\right]\,.
 \] 
We introduce
  the \emph{state operator} $\ZA  $ : $M^2_\zt\mapsto M^2_\zt$, the \emph{input operator} $\ZB $: $\zzr^m\mapsto M^2_\zt$  and the \emph{output operator } $\ZC $: $M^2_\zt\mapsto \zzr^p$:
 \[
 \begin{array}{l}
 \Dom\,\ZA  =\left \{ (\hat\xi,\tilde\xi)\in M^2_\zt\,,\quad \tilde \xi\in H^1(0,\zt)\,,\ \tilde \xi(0)=\hat\xi\right \}\\
\quad \ZA  \left[\begin{array}{c}
\hat\xi\\
\tilde \xi
\end{array}\right]
 =\left[\begin{array}{c}
A\hat\xi+\intzt N(s)\tilde\xi(s)\ZD s\\
-D_s\tilde \xi(s)
\end{array}\right] \,,\\
\quad
\ZB  u=\left[\begin{array}{c}
Bu\\
0
\end{array}\right] \qquad \ZC =\left[\begin{array}{cc}
C&0
\end{array}\right]\,.
\end{array}
 \]
  
  We search for a representation of the Riccati equation~(\ref{Eq:EquqRiccati})
  in terms of the state, input and outputs operators just introduced.
We verify that the following equality holds for every
$\ZOMq$ and $\Xi$ in $\Dom\,\ZA$:
\begin{multline}\ZLA{eq:RiccaOPERATform}  
 \dfrac{\ZD}{\ZD \zt} \ZL \ZOMq, \PP( \zt)\Xi\ZR+\ZL\ZA\ZOMq,\PP( \zt)\Xi\ZR+\ZL\PP( \zt)\ZOMq,\ZA \Xi\ZR\\
 -\ZL\ZB^*\PP( \zt)\ZOMq,B^*\PP( \zt)\Xi\ZR+\ZL \ZC \ZOMq,\ZC\Xi\ZR=0\qquad \PP(T)=0\,.
 \end{multline}
 
 \begin{Remark}{\rm
 We observe:
 \begin{enumerate}
 \item  the $L^2$ component of the elements of $\Dom\,\ZA$ belong to $H^1$ so that their derivative exist a.e. and are square integrable.
 \item the state, input and output operators do depend on the time $\zt$  since their definition depends on the space $M^2_\zt$, which depends on $\zt$. 
 
 For every $\zt$, the operator $\ZA $ is densely defined and closed. 
 
 \item different state space representations for systems with memory have been used. In particular, a representation via semigroup theory, see for example~\cite{Amendola2012,ENGELnagelLIBOsemigruppi}. The advantage of the semigroup approach is that it represents the system as a \emph{time invariante system} but this approach requires that the system be considered as a system with \emph{infinite memory:}   the memory is up to $-\ZIN$.
 The semigroup approach   is natural when studying the system on the entire time axis,
 as in the study of stability (see for example~\cite{Amendola2012} or~\cite[Chap.~5]{PandolfiLIBRO21}) but it is not natural when the system is considered on a finite time horizon, as in the study of controllability (see\cite[Chaps.~3-4]{PandolfiLIBRO21}) or of the regulator problem.  Even more so because we noted that the addition of the   memory of the configuration of the body   before the initial time does not add generality to the problem,
see Remark~\ref{remaSect1EFFETTOpassato}.

 \item the paper~\cite{acquistaBUCCI24} proposes an operator form of the Riccati equation which is different from~(\ref{eq:RiccaOPERATform}). Whether the  operators proposed in~\cite{acquistaBUCCI24}   have  any relation with  a state   representation of the system with memory is not investigated there.\zdia  
 \end{enumerate}
}\end{Remark}  

The validity of~(\ref{eq:RiccaOPERATform})
 is just a verification which uses the fact that the matrices $P_1(s,\zt)$ and $P_2(s,\nu,\zt)$ are differentiable in their variables.
 We confine ourselves to describe the key points.

First we compute the derivative of $ \zt\mapsto \ZL \ZOMq, \PP(\zt)\Xi\ZR$. When computing the derivatives of $\intzt P_1(\zt-s,\zt)\tilde\xi(s)\ZD s$ and of $\intzt P_2(\zt-s,\zt-\nu,\zt)\tilde\xi(\nu)\ZD\nu$ we use differentiability of the matrices $P_i$ respect to their entries and $\tilde\xi(0)=\hat\xi$, $\tilde\ZOM(0)=\hat\ZOM$, which is implied by the properties of $\Dom\ZA$.

We us $D_{,i}$ to denote the derivative respect to the $i-th$ variable while the derivative respect to $\zt$ of the function $P_0(\zt)$ (which depends on the sole variable $\zt$) is indicated with the superimposed dot. 

We insert the symbol  $(*)$  at the end of the lines which contain terms which do not appear in the form~(\ref{Eq:EquqRiccati}) of the Riccati equation found in the previous section.

We have:
\begin{equation}\ZLA{eq:ChStateRiccaDeriPt}
\begin{array}{ll}\displaystyle
\dfrac{\ZD}{\ZD \zt}\ZL\ZOM,\PP( \zt )\Xi\ZR\\
\displaystyle
= \ZL\hat\ZOM,\dot P_0(\zt)\hat\xi\ZR+\int_0^\zt \ZL \hat\ZOM,D_{,2}P_1(\zt-s,\zt)\tilde\xi(s)\ZD s \\
\displaystyle
+\ZL \hat\ZOM,P_1(0,\zt)\tilde\xi(\zt)\ZR+\int_0^\zt\ZL\hat\ZOM,D_{,1}P_1(\zt-s,\zt) \tilde\xi(s)\ZD s& {\bf (*)}\\
\displaystyle
+\ZL\tilde\ZOM(\zt),P_1^*(0,\zt)\hat\xi\ZR+\int_0^\zt\ZL\tilde\ZOM(s),D_{,1}P_1^*(\zt-s,\zt)\hat\xi\ZR\ZD s& {\bf (*)}\\
\displaystyle
+\int_0^\zt\ZL\tilde\ZOM(s),D_{,2}P_1^*(\zt-s,\zt)\hat\xi\ZR\ZD s
+\int_0^\zt\int_0^\zt\ZL\tilde\ZOM(s),D_{,3}P_2(\zt-s,\zt-\nu,\zt)\tilde\xi(\nu)\ZR\ZD\nu\,\ZD s\\
\displaystyle
+\int_0^\zt\ZL\tilde\ZOM(\zt),P_2(0,\zt-\nu,\zt)\tilde\xi(\nu)\ZR\ZD\nu+\int_0^\zt\ZL\tilde\ZOM(s),P_2(\zt-s,0,\zt)\tilde\xi(\zt)\ZR\ZD s& {\bf (*)}\\
\displaystyle
+\int_0^\zt\int_0^\zt \ZL\tilde\ZOM(s),[D_{,1}+D_{,2}]P_2(\zt-s,\zt-\nu,\zt)\tilde\xi(\nu)\ZR\ZD\nu\,\ZD s\,.& {\bf (*)}\\
\end{array}
\end{equation}

 Next we compute $\ZL \PP(\zt )\ZOMq,\ZA \Xi\ZR$. The derivative of $\tilde\xi$   in this expression  is absorbed by integrating by parts. We get:
 
\begin{multline}\ZLA{eq:ChStateRiccaprodAwithPt}
\begin{array}{ll}
\displaystyle 
\ZL \PP( \zt )\ZOMq,\ZA \Xi\ZR   &  \\
\displaystyle 
=\ZL \hat\ZOM , P_0( \zt )A\hat\xi\ZR+\intzt \ZL \hat\ZOM,P_0( \zt )N(s)\tilde\xi(s)\ZR\ZD s   &  \\
\displaystyle
+\intzt \ZL\tilde\ZOM(s),P_1^*( \zt -s, \zt )A\hat\xi\ZR\ZD s 
  +\intzt\intzt \ZL\tilde\ZOM(s),P_1^*( \zt -s, \zt )N(\nu)\tilde\xi(\nu)\ZR\ZD\nu\,\ZD s    & \\
  \displaystyle
  -\int_0^\zt \ZL P_2(0,\zt-s,\zt)\tilde\ZOM(s), \tilde\xi(\zt)\ZR\ZD s& (*)\\
  \displaystyle
 +\ZL \hat\ZOM,P_1( \zt , \zt )\hat\xi\ZR   +\intzt \ZL \tilde\ZOM(s),P_2^*( \zt  , \zt-s , \zt )\hat\xi\ZR\ZD s   &  \\
 \displaystyle 
-\intzt \ZL\hat \ZOM ,D_{,1} P_1( \zt -\nu, \zt )\tilde\xi(\nu)\ZD\nu
-\intzt\intzt\ZL\tilde\ZOM(s),D_{,1} P_2^*( \zt -\nu, \zt -s, \zt )\tilde\xi(\nu)\ZD\nu\,\ZD s\,.   & (*)
\end{array}
\end{multline}
We sum~(\ref{eq:ChStateRiccaprodAwithPt}) to~(\ref{eq:ChStateRiccaDeriPt}).
  The lines marked (*) in~(\ref{eq:ChStateRiccaprodAwithPt}) cancel part of the lines~(*) in~(\ref{eq:ChStateRiccaDeriPt}).  The  (*)-terms which remain inin~(\ref{eq:ChStateRiccaDeriPt}) are cancelled when we sum $\ZL\ZA\ZOMq, \PP( \zt )\Xi\ZR$.   The Riccati equation~(\ref{eq:RiccaOPERATform})  is verified by comparing with~(\ref{Eq:EquqRiccati}).

  Now we consider the contribution of the reference signal $y$ and we obtain the operator version of~(\ref{LeEquaSINTESdD}). For every $\zt\in[0,T]$ we introduce 
\[
\bd(\zt)=(\dd(\zt),\D(\zt-\cdot,\zt))\in M^2_\zt
\]
and we verify that the following  equality holds for every $\Xi\in\Dom\,\ZA$:
\begin{equation}
\ZLA{LeEquaSINTESdDopFORM}
\dfrac{\ZD}{\ZD\zt}\ZL\bd(\zt),\Xi\ZR=
 -\ZL\bd(\zt),\left (\ZA-\ZB\ZB^*\PP(\zt)\right )\Xi\ZR+\ZL y(\zt),\ZC\Xi\ZR\,,\quad \bd(T)=0\,.
\end{equation} 
In fact, by using~(\ref{LeEquaSINTESdD}) we compute that
\begin{align*}
\dfrac{\ZD}{\ZD\zt}\ZL\bd(\zt),\Xi\ZR=-\ZL \dd (\zt),A\hat\xi\ZR+\ZL \dd (\zt),BB^*P_0(\zt)\hat\xi\ZR+\ZL y(\zt),C\hat\xi\ZR \\
+\int_0^\zt  \ZL\D( \zt-s,\zt),(\mbox{\rm d}/\mbox{\rm d} s){\tilde\xi}(s)\ZD s 
-\int_0^\zt \ZL\dd(\zt),N(s)\tilde\xi(s)\ZR\ZD s\\
+\int_0^\zt\ZL \dd(\zt),BB^*P_1(\zt-s,s)\hat\xi\ZR\ZD s\,.
\end{align*}
The equality of this expression with the right side of~(\ref{LeEquaSINTESdDopFORM}) is easily verified.

\section{Conclusion} 
In this paper we extended the Fredholm equation approach to the study of the quadratic tracking problem on a finite time horizon for a linear systems with persistent memory and we obtained the synthesis of the optimal control.  The system is described by a Volterra integrodifferential equation. In the final section we introduced a new state space representation \emph{with finite memory} for the system and we recasted the synthesis in terms of a system of differential equations in this state space.


\begin{thebibliography}{10}

\bibitem{acquistaBUCCI24}
P.~Acquistapace and F.~Bucci.
\newblock Riccati-based solution to the optimal control of linear evolution
  equations with finite memory.
\newblock {\em Evol. Equ. Control Theory}, 13(1):26--66, 2024.

\bibitem{Amendola2012}
G.~Amendola, M.~Fabrizio, and J.~Golden.
\newblock {\em Thermodynamics of materials with memory: Theory and
  applications}.
\newblock {S}pringer, New York, 2012.

\bibitem{AthansFALB1966}
M.~Athans and P.~L. Falb.
\newblock {\em Optimal control. {A}n introduction to the theory and its
  applications}.
\newblock McGraw-Hill Book Co., New York-Toronto-London, 1966.

\bibitem{DAVISfoundatCONTROL02}
J.~H. Davis.
\newblock {\em Foundations of deterministic and stochastic control}.
\newblock Systems \& Control: Foundations \& Applications. Birkh\"{a}user
  Boston, Inc., Boston, MA, 2002.

\bibitem{ENGELnagelLIBOsemigruppi}
K.-J. Engel and R.~Nagel.
\newblock {\em One-parameter semigroups for linear evolution equations}, volume
  194 of {\em Graduate Texts in Mathematics}.
\newblock Springer-Verlag, New York, 2000.
\newblock With contributions by S. Brendle, M. Campiti, T. Hahn, G. Metafune,
  G. Nickel, D. Pallara, C. Perazzoli, A. Rhandi, S. Romanelli and R.
  Schnaubelt.

\bibitem{FagnaniMaksPandolfi}
F.~Fagnani, V.~Maksimov, and L.~Pandolfi.
\newblock A recursive deconvolution approach to disturbance reduction.
\newblock {\em IEEE Trans. Automat. Control}, 49(6):907--921, 2004.

\bibitem{FagnaniPandJIIPP03}
F.~Fagnani and L.~Pandolfi.
\newblock On the solution of a class of {V}olterra integral equations of the
  first kind.
\newblock {\em J. Inverse Ill-Posed Probl.}, 11(5):485--503, 2003.

\bibitem{MFAGnanIPandoLfiIPvoltyeDecon03}
F.~Fagnani and L.~Pandolfi.
\newblock A recursive algorithm for the approximate solution of {V}olterra
  integral equations of the first kind of convolution type.
\newblock {\em Inverse Problems}, 19(1):23--47, 2003.

\bibitem{ShuoLin23}
S.~Han, P.~Lin, and J.~Yong.
\newblock Causal state feedback representation for linear quadratic optimal
  control problems of singular {V}olterra integral equations.
\newblock {\em Math. Control Relat. Fields}, 13(4):1282--1317, 2023.

\bibitem{MaksimovLIBRO02}
V.~I. Maksimov.
\newblock {\em Dynamical inverse problems of distributed systems}.
\newblock Inverse and Ill-posed Problems Series. VSP, Utrecht, 2002.

\bibitem{ManitiusOPTIMALtrieste76}
A.~Manitius.
\newblock Optimal control of hereditary systems.
\newblock In {\em Control theory and topics in functional analysis ({I}nternat.
  {S}em., {I}nternat. {C}entre {T}heoret. {P}hys., {T}rieste, 1974), {V}ol.
  {III}}, pages 43--178. Internat. Atomic Energy Agency, Vienna, 1976.

\bibitem{PandolfiINTjCONTreject}
L.~Pandolfi.
\newblock Adaptive recursive deconvolution and adaptive noise cancellation.
\newblock {\em Internat. J. Control}, 80(3):403--415, 2007.

\bibitem{PandolfiVOLTERRAieee18}
L.~Pandolfi.
\newblock The quadratic regulator problem and the {R}iccati equation for a
  process governed by a linear {V}olterra integrodifferential equations.
\newblock {\em IEEE Trans. Automat. Control}, 63(5):1517--1522, 2018.

\bibitem{PandolfiLIBRO21}
L.~Pandolfi.
\newblock {\em Systems with persistent memory---controllability, stability,
  identification}, volume~54 of {\em Interdisciplinary Applied Mathematics}.
\newblock Springer, Cham, [2021] \copyright 2021.

\bibitem{PandolfiPRIOLAtracking}
L.~Pandolfi and E.~Priola.
\newblock Tracking control of parabolic systems.
\newblock In {\em System modeling and optimization}, volume 166 of {\em IFIP
  Int. Fed. Inf. Process.}, pages 135--146. Kluwer Acad. Publ., Boston, MA,
  2005.

\bibitem{pritchardLUO96}
A.~J. Pritchard and Y.~You.
\newblock Causal feedback optimal control for {V}olterra integral equations.
\newblock {\em SIAM J. Control Optim.}, 34(6):1874--1890, 1996.

\bibitem{PruessLIBRO1993}
J.~Pr\"uss.
\newblock {\em Evolutionary Integral Equations and Applications}.
\newblock Birkh\"auser, Basel, 1993.

\bibitem{SontagLIBRO98}
E.~D. Sontag.
\newblock {\em Mathematical control theory}, volume~6 of {\em Texts in Applied
  Mathematics}.
\newblock Springer-Verlag, New York, second edition, 1998.
\newblock Deterministic finite-dimensional systems.

\bibitem{WANG18ESAIm}
T.~Wang.
\newblock Linear quadratic control problems of stochastic {V}olterra integral
  equations.
\newblock {\em ESAIM Control Optim. Calc. Var.}, 24(4):1849--1879, 2018.

\end{thebibliography}
\enddocument